\documentclass[11pt, eqno]{article}

\usepackage{amssymb,latexsym}
\usepackage{amsmath,latexsym}
\usepackage[top=1in,bottom=1in,left=1.3in,right=1.3in]{geometry}
\usepackage{amscd}
\usepackage{amsthm}

\newtheorem{theorem}{Theorem}[section]
\newtheorem{cor}[theorem]{Corollary}

\newtheorem{thm}[theorem]{Theorem}

\newtheorem{rem}[theorem]{Remark}
\newtheorem{lemma}[theorem]{Lemma}

\newtheorem{prop}[theorem]{Proposition}

\begin{document}

\title{Real hypersurfaces of complex quadric in terms of star-Ricci tensor}

\author{\\  Xiaomin Chen
\thanks{
The author is supported by the Science Foundation of China
University of Petroleum-Beijing(No.2462015YQ0604) and partially supported
by  the Personnel Training and Academic
Development Fund (2462015QZDX02).
 }\\
{\normalsize College of  Science, China University of Petroleum (Beijing),}\\
{\normalsize Beijing, 102249, China}\\
{\normalsize xmchen@cup.edu.cn}}
\maketitle \vspace{-0.1in}


\abstract{In this article, we introduce the notion of star-Ricci tensors in the real hypersurfaces of complex quadric $Q^m$.
It is proved that there exist no Hopf hypersurfaces in $Q^m,m\geq3$,
with commuting star-Ricci tensor or parallel star-Ricci tensor. As a
generalization of star-Einstein metric, star-Ricci solitons on $M$ are considered.
In this case we show that $M$ is an open part
of a tube around a totally geodesic $\mathbb{C}P^\frac{m}{2}\subset Q^{m},m\geq4$.}
 \vspace{-0.1in}
\medskip\vspace{12mm}

\noindent{\it Keywords}:  Hopf hypersurface; complex quadric; commuting star-Ricci tensor; parallel star-Ricci tensor; star-Ricci soliton
  \vspace{2mm}

\noindent{\it MSC}: 53C40, 53C15 \vspace{2mm}

\section{Introduction}

The complex quadric $Q^m$ is a Hermitian symmetric space $\mathrm{SO}_{m+2}/\mathrm{SO}_m\mathrm{SO}_2$
with rank two in the class of compact type. It can be regarded as a complex
hypersurface of complex projective space $\mathbb{C}P^{m+1}$. Also, the
complex quadric $Q^m$ can be regarded as a kind of real Grassmannian
manifolds of compact type with rank two. In the complex quadric $Q^m$
there are two important geometric structures, a complex conjugation
structure $A$ and K\"ahler structure $J$, with each other being
anti-commuting, that is, $AJ=-JA$. Another distinguish geometric structure in $Q^m$ is a parallel rank two
vector bundle $\mathfrak{U}$ which contains an $S^1$-bundle of
real structures, that is, complex conjugations $A$ on the tangent spaces of
$Q^m$. Here the parallel vector bundle $\mathfrak{U}$ means
that $(\widetilde{\nabla}_XA)Y=q(X)AY$ for all $X,Y\in T_zQ^m,z\in Q^m,$
where $\widetilde{\nabla}$ and $q$ denote a connection and a certain 1-form on $T_zQ^m,$ respectively.

Recall that a nonzero tangent vector $W\in T_zQ^m,z\in Q^m$, is called \emph{singular}
 if it is tangent to more than one maximal flat in $Q^m$. There are two types of
 singular tangent vectors for the complex quadric $Q^m$:
\begin{enumerate}
  \item If there exists a conjugation $A\in\mathfrak{U}$ such that $W\in V(A)$,
  then $W$ is singular. Such a singular tangent vector is called $\mathfrak{U}$-principal.
  \item If there exist a conjugation $A\in\mathfrak{U}$ and orthonormal
  vectors $X,Y\in V(A)$ such that $W/\|W\| =(X+JY)/\sqrt{2},$ then $W$ is
  singular. Such a singular tangent vector is called $\mathfrak{U}$-isotropic.
\end{enumerate}

Let $M$ be a real hypersurface of $Q^m$. The K\"ahler structure $J$ on $Q^m$
induces a structure vector field $\xi$ called {\it Reeb vector field} on
$M$ by $\xi:=-JN$, where $N$ is a local unit normal vector field of $M$ in $Q^m$.
It is well-known that there is an almost contact structure $(\phi,\eta,\xi,g)$
 on $M$ induced from complex quadric. Moreover, if the Reeb vector field $\xi$
  is invariant under the shape operator $S$, i.e. $S\xi=\alpha\xi$, where
  $\alpha=g(S\xi,\xi)$ is a smooth function, then $M$ is
said to be a \emph{Hopf hypersurface}. For the real Hopf hypersurfaces of
 complex quadric many characterizations were obtained by Suh
  (see \cite{Suh1,Suh2,Suh3,Suh4,Suh5} etc.). In particular, we note that
  Suh in \cite{Suh1} introduced parallel Ricci tensor, i.e.$\nabla \mathrm{Ric}=0$,
   for the real hypersurfaces in $Q^m$ and gave a complete
   classification for this case.
 In addition, if the real hypersurface $M$ admits commuting Ricci
 tensor, i.e. $\mathrm{Ric}\circ\phi=\phi\circ \mathrm{Ric}$, Suh also proved the followings:
\begin{thm}[\cite{Suh5}] Let M be a real hypersurface of the complex
quadric $Q^m,m\geq3$, with commuting Ricci tensor. Then the unit
normal vector field $N$ of $M$ is either $\mathfrak{U}$-principal or $\mathfrak{U}$-isotropic.
\end{thm}
\begin{thm}[\cite{Suh5}]
 There exist no Hopf real hypersurfaces in the complex quadric $Q^m, m\geq4$, with commuting and
parallel Ricci tensor.
\end{thm}

Since the Ricci tensor of an Einstein hypersurface in the complex quadric
$Q^m$ is a constant multiple of $g$, it  satisfies naturally commuting and
parallelism. Thus we have the following.
\begin{cor}[\cite{Suh5}]\label{C1}
 There exist no Hopf Einstein real hypersurfaces in the complex quadric $Q^m,m\geq4.$
\end{cor}

As a generalization of Einstein metrics, recently Suh in \cite{Suh6}  has shown
 a complete classification of Hopf hypersurfaces with a Ricci soliton, which is given by
\begin{equation*}
\frac{1}{2}(\mathfrak{L}_Wg)(X,Y)+\mathrm{Ric}(X,Y)=\lambda g(X,Y).
\end{equation*}
Here $\lambda$ is a constant and $W$ is a vector field on $M$, which are said
to be  \emph{Ricci soliton constant}  and \emph{potential vector field}, respectively,
and $\mathfrak{L}_W$ denotes the Lie derivative along the direction of the vector field $W$.

Notice that, as the corresponding of Ricci tensor, Tachibana \cite{T}
introduced the idea of star-Ricci tensor. These ideas
 apply to almost contact metric manifolds, and in particular,
 to real hypersurfaces in complex space forms by Hamada in \cite{HA}.
The star-Ricci tensor $\mathrm{Ric}^*$  is defined by
\begin{equation}\label{1}
  \mathrm{Ric}^*(X,Y)=\frac{1}{2}\mathrm{trace}\{\phi\circ R(X,\phi Y)\},\quad\text{for all}\, X,Y\in TM.
\end{equation}
If the star-Ricci tensor is a constant multiple of $g(X,Y)$ for all $X,Y$ orthogonal
to $\xi$, then $M$ is said to be a \emph{star-Einstein manifold}.
 Hamada gave a classification of star-Einstein hypersurfaces of
 $\mathbb{C}P^n$ and $\mathbb{C}H^n$, and further Ivey
 and Ryan updated and refined the work of Hamada in 2011(\cite{IR}).

Motivated by the present work, in this paper we introduce the notion of star-Ricci tensor
in the real hypersurfaces of complex quadric $Q^m$ and study the characterizations of a
real Hopf hypersurface whose star-Ricci tensor satifies certain conditions.

First we consider
the real hypersurface with commuting star-Ricci tensor,
i.e. $\phi\circ \mathrm{Ric}^*=\mathrm{Ric}^*\circ\phi$.  We assert the following:
\begin{thm}\label{T1*}
There exist no Hopf hypersurfaces of $Q^m,m\geq3$, with commuting star-Ricci tensor.
\end{thm}

For the Hopf hypersurfaces of $Q^m,m\geq3,$ with parallel
star-Ricci tensor, we also prove the following non-existence.
\begin{thm}\label{T3}
There exist no Hopf hypersurfaces of $Q^m,m\geq3$, with parallel star-Ricci tensor.
\end{thm}

As the generalization of star-Einstein metric Kaimakamis and Panagiotidou \cite{KP}
 introduced a so-called \emph{star-Ricci soliton}, that is, a
Riemannain metric $g$ on $M$ satisfying
\begin{equation}\label{2}
 \frac{1}{2}\mathcal{L}_W g+\mathrm{Ric}^*=\lambda g.
\end{equation}
In this case we obtain the following characterization:
\begin{thm}\label{T2}
Let $M$ be a real hypersurface in $Q^m,m\geq4$, admitting a star-Ricci
soliton with potential vector field $\xi$, then $M$ is an open part
of a tube around a totally geodesic $\mathbb{C}P^\frac{m}{2}\subset Q^{m}$.
\end{thm}

This paper is organized as follows. In Section 2 and Section 3, some
basic concepts and formulas for real hypersurfaces in complex
quadric are presented. In Section 4 we consider Hopf hypersurfaces
 with commuting star-Ricci tensor and give the proof of Theorem \ref{T1*}.
In Section 5 we will prove Theorem \ref{T3}. At last we assume that a
 Hopf hypersurface admits star-Ricci soliton and give the proof of Theorem \ref{T2} as the Section 6.

\section{The complex quadric }

 In this section we will summarize some basic notations and formulas about
 the complex quadric $Q^m$.
 For more detail see \cite{BS,BS2,KN,K}.
The complex quadric $Q^m$ is the hypersurface of complex projective space $\mathbb{C}P^{m+1}$,
which is defined by $z_1^2+\cdots+z_{m+2}^2=0$, where $z_1,\cdots,z_{m+2}$ are homogeneous
coordinates on $\mathbb{C}P^{m+1}.$ In the complex quadric it
is equipped with a Riemannian metric $\widetilde{g}$ induced from
the Fubini-Study metric on $\mathbb{C}P^{m+1}$ with constant holomorphic
sectional curvature 4. Also the K\"ahler structure on $\mathbb{C}P^{m+1}$ induces
canonically a K\"ahler structure $(J,\widetilde{g})$ on the complex quadric $Q^m.$

The complex projective space $\mathbb{C}P^{m+1}$ is a Hermitian symmetric
space of the special unitary group $\mathrm{SU}_{m+2}$, i.e. $\mathbb{C}P^{m+1}=\mathrm{SU}_{m+2}/S(U_1U_{m+1})$.
The special orthogonal group $\mathrm{SO}_{m+2}\subset \mathrm{SU}_{m+2}$ acts on
$\mathbb{C}P^{m+1}$ with cohomogeneity one. The orbit containing $o$ is a totally geodesic real projective space
$\mathbb{R}P^{m+1}\subset\mathbb{C}P^{m+1}$, where
$o=[0,\cdots,0,1]\in \mathbb{C}P^{m+1}$ is the fixed point of the action
 of the stabilizer $S(U_{m+1}U_1)$. We can identify $Q^m$ with a homogeneous space
$\mathrm{SO}(m+2)/\mathrm{SO}_2\mathrm{SO}_m$, which is the second singular orbit of this action.
Such a homogeneous space model leads to the geometric interpretation
of the complex quadric $Q^m$ as the Grassmann manifold $G^+_2(\mathbb{R}^{m+2})$
of oriented 2-planes in $\mathbb{R}^{m+2}$. From now on we always assume $m\geq3$
because it is well known that $Q^1$ is isometric
to a sphere $S^2$ with constant curvature and $Q^2$ is isometric to the Riemannian
product of two 2-spheres with constant curvature.

For a unit normal vector $\rho$ of $Q^m$ at a point $z\in Q^m$ we denote by $A=A_\rho$
the shape operator of $Q^m$ in $\mathbb{C}P^{m+1}$ with respect to $\rho$, which is an
involution on the tangent space $T_zQ^m$, and the tangent space can be decomposed as
\begin{equation*}
T_zQ^m=V(A_\rho)\oplus JV(A_\rho),
\end{equation*}
where $V(A_\rho)$ is the $(+1)$-eigenspace and $JV(A_\rho)$ is the $(-1)$-eigenspace of $A_\rho$.
This means that the shape operator $A$ defines a real structure
on $T_zQ^m,$ equivalently, $A$ is a complex conjugation.
Since the real codimension of $Q^m$ in $\mathbb{C}P^{m+1}$ is 2,
this induces an $S^1$-subbundle $\mathfrak{U}$
of the endomorphism bundle $\mathrm{End}(TQ^m)$ consisting of complex conjugations.
Notice that $J$ and each complex conjugation $A\in\mathfrak{U}$ anti-commute, i.e. $AJ=-JA.$

\section{Real hypersurface of complex quadric and its star-Ricci tensor}

Let $M$ be an immersed real hypersurface of $Q^m$ with induced metric $g$. There exists a local defined
unit normal vector field $N$ on $M$ and we write $\xi:=- JN$ by the structure vector field of $M$.
 An induced one-form $\eta$ is defined by
$\eta(\cdot)=\widetilde{g}(J\cdot,N)$, which is dual to $\xi$.  For any vector field $X$ on $M$ the tangent part of $JX$
is denoted by $\phi X=JX-\eta(X)N$. Moreover, the following identities hold:
\begin{align}
&\phi^2=-Id+\eta\otimes\xi,\quad\eta\circ \phi=0,\quad\phi\circ\xi=0,\quad\eta(\xi)=1,\label{3}\\
&g(\phi X,\phi Y)=g(X,Y)-\eta(X)\eta(Y),\quad g(X,\xi)=\eta(X),\label{4}
\end{align}
where $X,Y\in\mathfrak{X}(M)$. By these formulas, we know that $(\phi,\eta,\xi,g)$ is an almost
contact metric structure on $M$. The tangent bundle $TM$ can be decomposed
as $TM=\mathcal{C}\otimes\mathbb{R}\xi$, where $\mathcal{C}=\ker\eta$ is
 the maximal complex subbundle of $TM$. Denote by $\nabla, S$ the
 induced Riemannian connection and the shape operator on $M$, respectively.
Then the Gauss and Weingarten formulas are given respectively by
\begin{equation}\label{5}
\widetilde{\nabla}_XY=\nabla_XY+g(SX,Y)N,\quad\widetilde{\nabla}_XN=-SX,
\end{equation}
 where $\widetilde{\nabla}$ is the connection on $Q^m$ with respect to $\widetilde{g}$.
Also, we have
\begin{equation}\label{6}
  (\nabla_X\phi)Y=\eta(Y)SX-g(SX,Y)\xi,\quad\nabla_X\xi=\phi SX.
\end{equation}

The curvature tensor $R$ and Codazzi equation of $M$ are given respectively as follows(see \cite{Suh1}):
 \begin{align}\label{7}
R(X,Y)Z=&g(Y,Z)X-g(X,Z)Y+g(\phi Y,Z)\phi X-g(\phi X,Z)\phi Y-2g(\phi X,Y)\phi Z\nonumber\\
&+g(AY,Z)AX-g(AX,Z)AY+g(JAY,Z)JAX-g(JAX,Z)JAY\nonumber\\
&+g(SY,Z)SX-g(SX,Z)SY,
\end{align}
\begin{equation}\label{8}
\begin{aligned}
 g((\nabla_XS)Y-(\nabla_YS)X,Z)=&\eta(X)g(\phi Y,Z)-\eta(Y)g(\phi X,Z)-2\eta(Z)g(\phi X,Y)\\
&+g(X,AN)g(AY,Z)-g(Y,AN)g(AX,Z)\\
&+g(X,A\xi)g(JAY,Z)-g(Y,A\xi)g(JAX,Z)
 \end{aligned}
 \end{equation}
for any vector fields $X,Y,Z$ on $M$.

At each point $z\in M$ we denote
\begin{equation*}
  \mathcal{Q}_z=\{X\in T_zM|\;AX\in T_zM\;\text{for all}\; A\in\mathfrak{U}_z\}
\end{equation*}
by a maximal $\mathfrak{U}$-invariant subspace of $T_zM$. For the subspace the following lemma was proved.
\begin{lemma}[see \cite{Suh2}] For each $z\in M$ we have
\begin{itemize}
  \item If $N_z$ is $\mathfrak{U}$-principal, then $\mathcal{Q}_z=\mathcal{C}_z$.
  \item If $N_z$  is not $\mathfrak{U}$-principal, there exist a
  conjugation $A\in\mathfrak{U}$ and orthonormal vectors $X,Y\in V(A)$
  such that $N_z=\cos(t)X+\sin(t)JY$ for some $t\in(0,\frac{\pi}{4}]$.
  Then we have $\mathcal{Q}_z=\mathcal{C}_z\ominus\mathbb{C}(JX+Y)$.
\end{itemize}
\end{lemma}

For each point $z\in M$ we choose $A\in\mathfrak{U}_z$, then
there exist two orthonormal vectors $Z_1,Z_2\in V(A)$
such that
\begin{equation}\label{9}
\left\{
\begin{array}{ll}
N&=\cos(t)Z_1+\sin(t)JZ_2,\\
AN&=\cos(t)Z_1-\sin(t)JZ_2,\\
\xi&=\sin(t)Z_2-\cos(t)JZ_1,\\
A\xi&=\sin(t)Z_2+\cos(t)JZ_1
  \end{array}
\right.
\end{equation}
for $0\leq t\leq \frac{\pi}{4}$ (see \cite[Proposition 3]{R}). From this we get $g(AN,\xi)=0.$

In the real hypersurface $M$ we introduce the star-Ricci tensor $\mathrm{Ric}^*$ defined by
\begin{equation*}
\mathrm{Ric}^*(X,Y)=\frac{1}{2}\mathrm{trace}\{\phi\circ R(X,\phi Y)\},\quad\hbox{for all}\, X,Y\in TM.
\end{equation*}
Taking a local frame $\{e_i\}$ of $M$ such that $e_1=\xi$ and using \eqref{4}, we derive from \eqref{7}
\begin{align*}
 &\sum_{i=1}^{2m-1}g(\phi\circ R(X,\phi Y)e_i,e_i)\\
 =&g(\phi Y,\phi X)-g(X,\phi^2 Y)+g(\phi^2 Y,\phi^2 X)-g(\phi X,\phi^3 Y)+2(2m-2)g(\phi X,\phi Y)\\
&+g(A\phi Y,\phi AX)-g(AX,\phi A\phi Y)+g(JA\phi Y,\phi JAX)-g(JAX,\phi JA\phi Y)\\
&+g(S\phi Y,\phi SX)-g(SX,\phi S\phi Y)\\
   =&4mg(\phi X,\phi Y) -2g(AX,\phi A\phi Y)+2g(JA\phi Y,\phi JAX)-2g(SX,\phi S\phi Y).
\end{align*}
In view of \eqref{1}, the star-Ricci tensor is given by
\begin{equation}\label{10}
\begin{aligned}
\mathrm{Ric}^*(X,Y)=&2mg(\phi X,\phi Y) -g(AX,\phi A\phi Y)\\
&+g(JA\phi Y,\phi JAX)-g(SX,\phi S\phi Y).
\end{aligned}
\end{equation}
Since $AJ=-JA$ and $\xi=-JN$, we have
\begin{align*}
  JA\phi Y=&-AJ\phi Y=AY-\eta(Y)A\xi, \\
  \phi JAX=&J(JAX)-\eta(JAX)N=-AX+g(N,AX)N.
\end{align*}
Then
\begin{align}\label{11}
  g(JA\phi Y,\phi JAX)=&-g(AX,AY)+\eta(Y)\eta(X)+g(N,AX)g(AY,N)\nonumber\\
  =&g(\phi^2X,Y)+g(N,AX)g(AY,N).
\end{align}
Because \begin{align*}
          JA\phi Y=&\phi A\phi Y+\eta(A\phi Y)N \\
          = & \phi A\phi Y+g(\xi,AJY-\eta(Y)AN)N\\
          =&\phi A\phi Y+g(J\xi,AY)N\\
          =&\phi A\phi Y+g(N,AY)N,
        \end{align*}
 we have
\begin{align}\label{12}
  g(AX,\phi A\phi Y) & =g(AX,JA\phi Y-g(N,AY)N) \nonumber\\
                 & =g(AX,JA\phi Y)-g(N,AY)g(AX,N)\nonumber\\
                 &=-g(\phi^2X,Y)-g(N,AY)g(AX,N).
\end{align}
Thus substituting \eqref{11} and \eqref{12} into \eqref{10} implies
\begin{align}\label{13}
  \mathrm{Ric}^*(X,Y) =&-2(m-1)g(\phi^2 X,Y)-2g(N,AX)g(AY,N)-g((\phi S)^2X,Y)
\end{align}
for all $X,Y\in TM$.

In the following we always assume that $M$ is a Hopf hypersurface in $Q^m$,
i.e. $S\xi=\alpha\xi$ for a smooth function $\alpha=g(S\xi,\xi).$
As in \cite{Suh1}, since $g(AN,\xi)=0$, by taking $Z=\xi$
in the Codazzi equation \eqref{8},  we have
\begin{align*}
&g((\nabla_XS)Y-(\nabla_YS)X,\xi)\\
=&-2g(\phi X,Y)+2g(X,AN)g(AY,\xi)-2g(Y,AN)g(AX,\xi).
\end{align*}

On the other hand,
\begin{align*}
&g((\nabla_XS)Y-(\nabla_YS)X,\xi)\\
=& g((\nabla_XS)\xi,Y)-g((\nabla_YS)\xi,X)\\
=&(X\alpha)\eta(Y)-(Y\alpha)\eta(X)+\alpha g((\phi S+S\phi)X,Y)-2g(S\phi SX,Y).
\end{align*}
Comparing the previous two equations and putting $X=\xi$ gives
\begin{equation}\label{14}
  Y\alpha=(\xi\alpha)\eta(Y)+2g(Y,AN)g(\xi,A\xi).
\end{equation}
Moreover, we have the following.
\begin{lemma}(\cite[Lemma 4.2]{Suh2})
 Let $M$ be a Hopf hypersurface in $Q^m$ with (local) unit normal vector
field $N$. For each point in $z\in M$ we choose $A\in\mathfrak{U}_z$ such that $N_z=\cos(t)Z_1+
\sin(t)JZ_2$ holds for some orthonormal vectors $Z_1,Z_2\in V(A)$ and $0\leq t\leq \frac{\pi}{4}$. Then
\begin{equation}\label{15}
\begin{aligned}
0=&2g(S\phi SX,Y)-\alpha g((\phi S+S\phi )X,Y)-2g(\phi X,Y)\\
&+2g(X,AN)g(Y,A\xi)-2g(Y,AN)g(X,A\xi)\\
&+2g(\xi,A\xi)\{g(Y,AN)\eta(X)-g(X,AN)\eta(Y)\}
 \end{aligned}\end{equation}
holds for all vector fields $X,Y$ on $M$.
\end{lemma}
From this lemma we can prove the following.
\begin{lemma}\label{L2.1}
Let $M$ be a Hopf hypersurface in complex quadric $Q^m$, then
\begin{equation}\label{16}
  (\phi S)^2=(S\phi)^2.
\end{equation}
\end{lemma}
\proof From the equation \eqref{15} we assert the followings:
\begin{align}
g((S\phi)^2X,Y)=&\frac{1}{2}\alpha g((\phi S+S\phi )\phi X,Y)+g(\phi^2X,Y)-g(\phi X,AN)g(Y,A\xi)\nonumber\\
&+g(\phi X,A\xi)g(Y,AN)+g(\xi,A\xi)g(\phi X,AN)\eta(Y),\nonumber\\
g((\phi S)^2X,Y)=&\frac{1}{2}\alpha g(\phi(\phi S+S\phi )X,Y)+g(\phi^2 X,Y)-g(X,AN)g(\phi A\xi,Y)\nonumber\\
&+g(X,A\xi)g(\phi AN,Y)-g(\xi,A\xi)\eta(X)g(\phi AN,Y).\label{17}
 \end{align}
 Thus  we obtain
 \begin{align*}
   g((S\phi)^2X-(\phi S)^2X,Y)= &-g(\phi X,AN)g(Y,A\xi)+g(\phi X,A\xi)g(Y,AN)\nonumber\\
&+g(\xi,A\xi)g(\phi X,AN)\eta(Y)+g(X,AN)g(\phi A\xi,Y)\nonumber\\
&-g(X,A\xi)g(\phi AN,Y)+g(\xi,A\xi)\eta(X)g(\phi AN,Y)\\
=&\eta(X)g(AN,N)g(Y,A\xi)-g(\xi,A\xi)g(X,A\xi)\eta(Y)\nonumber\\
&-g(X,A\xi)\eta(Y)g(AN,N)+g(\xi,A\xi)\eta(X)g(Y,A\xi)\\
= &\Big(\eta(X)g(A\xi,Y)-g(X,A\xi)\eta(Y)\Big)\Big(g(AN,N)+g(\xi,A\xi)\Big).
 \end{align*}
 Here we have used the following relations:
 \begin{align}
   g(A\xi,\phi X) & = g(A\xi,JX-\eta(X)N)=g(AN,X), \label{18}\\
   g(A\phi X,N) & =g(AJX-\eta(X)AN,N)=-g(X,A\xi)-\eta(X)g(AN,N)\label{19}.
 \end{align}
From \eqref{9}, we get $g(AN,N)+g(\xi,A\xi)=0$, which yields \eqref{16}.\qed

 \section{Proof of Theorem \ref{T1*} }

  In this section we suppose that $M$ is a real Hopf hypersurface  with
  commuting star-Ricci tensor, that is, $\phi\circ \mathrm{Ric}^*=\mathrm{Ric}^*\circ\phi$.
  Making use of \eqref{13}, a straightforward computation gives
\begin{align*}
 0=& g((\phi\circ \mathrm{Ric}^*-\mathrm{Ric}^*\circ\phi)X,Y)\\
=&-\mathrm{Ric}^*(X,\phi Y)-\mathrm{Ric}^*(\phi X,Y)\\
  =&2g(N,AX)g(A\phi Y,N)+2g(N,A\phi X)g(AY,N)\nonumber\\
  &+g(\phi[(S\phi)^2-(\phi S)^2]X,Y).
\end{align*}
Thus Lemma \ref{L2.1} implies
 \begin{align*}
g(N,AX)g(A\phi Y,N)+g(N,A\phi X)g(AY,N)=0.
 \end{align*}
Replacing $X$ and $Y$ by $\phi X$ and $\phi Y$ respectively gives
 \begin{align*}
  g(N,A\phi X)g(Y,AN)+g(X,AN)g(A\phi Y,N)=0.
 \end{align*}
Now, if $X=Y$, we find $g(AN,\phi X)g(AN,X)=0$ for all vector field $X$ on $M$,
which means $AN=N$. Therefore we prove the following.
\begin{lemma}\label{L1}
Let $M$ be a Hopf hypersurface of complex quadric $Q^m, m\geq3$, with commuting
star-Ricci tensor. Then the unit normal vector field $N$ is $\mathfrak{U}$-principal.
\end{lemma}

 In terms of \eqref{17}, the star-Ricci tensor \eqref{13} becomes
 \begin{align*}
  \mathrm{Ric}^*(X,Y) =&(-2m+1)g(\phi^2 X,Y)-\frac{1}{2}\alpha g(\phi(\phi S+S\phi )X,Y).
\end{align*}
Moreover, from \eqref{15} we obtain
 \begin{align*}
  \mathrm{Ric}^*(X) =&(-2m+1)\phi^2 X-\frac{1}{2}\alpha\phi(\phi S+S\phi )X\\
  =&(-2m+1)\phi^2 X-\frac{1}{2}\alpha\phi^2SX-\frac{1}{4}\alpha^2 (\phi S+S\phi )X-\frac{1}{2}\alpha\phi X.
\end{align*}
By virtue of \cite[Lemma 4.3]{Suh1} and Lemma \ref{L1}, it implies that
$\alpha$ is constant. If $\alpha\neq0$, making use of the previous formula, we conclude that
\begin{equation*}
0=\phi \mathrm{Ric}^*(X)-\mathrm{Ric}^*(\phi X)= \frac{1}{2}\alpha(\phi SX-S\phi X)
\end{equation*}
for all $X\in TM$. That means that the Reeb flow is isometric. In view
of \cite[Proposition 6.1]{BS2}, the normal vector field $N$ is
isotropic everywhere, which is contradictory with Lemma \ref{L1}.
Hence $\alpha=0$ and the star-Ricci tensor becomes
 \begin{align}\label{20}
  \mathrm{Ric}^*(X,Y) =&(-2m+1)g(\phi^2 X,Y).
\end{align}

Now replacing $X$ and $Y$ by $\phi X$ and $\phi Y$ respectively in \eqref{13} and using \eqref{20}, we get
 \begin{align*}
  (2m-1)(\phi X,\phi Y) =&2(m-1)g(X,\phi Y)-2g(N,A\phi X)g(A\phi Y,N)-g((S\phi)^2 X,Y).
\end{align*}
Interchanging $X$ and $Y$ and applying the resulting equation to subtract the pervious equation, we obtain
 \begin{align*}
  g((S\phi)^2 X-(\phi S)^2X,Y)=&4(m-1)g(X,\phi Y).
\end{align*}
So from Lemma \ref{L2.1}, we conclude that
\begin{align*}
 4(m-1)g(X,\phi Y)=0,
 \end{align*}
which is impossible since $m\geq3$. We finish the proof of Theorem \ref{T1*}.\qed
\begin{rem}
Formula \eqref{20} with $X,Y\in\mathcal{C}$, we have $\mathrm{Ric}^*(X,Y)=(2m-1)g(X,Y)$,
namely $M$ is star-Einstein, thus we have proved that there exist no
star-Einstein Hopf hyersurfaces in complex quadric $Q^m,m\geq3$,
which is analogous to the conclusion of Corallary \ref{C1} in the introduction.
\end{rem}

\section{Proofs of Theorem \ref{T3}}
 In this section we assume $M$ is a Hopf hypersurface of $Q^m,m\geq3$,
 with parallel star-Ricci tensor. In order to prove Theorem \ref{T3},
 we first prove the following lemma.
  \begin{lemma}
 Let $M$ be a Hopf hypersurface of $Q^m,m\geq3,$ with parallel star-Ricci
 tensor. Then the unit normal vector $N$ is either $\mathfrak{U}$-principal or $\mathfrak{U}$-isotropic.
 \end{lemma}
 \proof Since $\nabla \mathrm{Ric}^*=0$, differentiating equation \eqref{13} covariantly along vector field $Z$ gives
\begin{align*}
  0=&2(m-1)g((\nabla_Z\phi)\phi X+\phi(\nabla_Z\phi)X,Y)\\
&+2g(\widetilde{\nabla}_ZN,AX)g(AY,N)+2g(N,(\widetilde{\nabla}_ZA)X)g(AY,N)\\
&+2g(\widetilde{\nabla}_ZN,AY)g(AX,N)+2g(N,(\widetilde{\nabla}_ZA)Y)g(AX,N)\\
&+g((\nabla_Z\phi) S\phi SX,Y)+g(\phi(\nabla_ZS)\phi SX,Y)\\
&+g(\phi S(\nabla_Z\phi) SX,Y)+g(\phi S\phi(\nabla_ZS)X,Y).
\end{align*}
Here we have used $(\widetilde{\nabla}_ZA)X=q(Z)AX$ for a certain 1-form $q$ as in the introduction.
Moreover, by \eqref{5} we have
\begin{align}\label{21}
  0=&-2(m-1)g(SZ,\phi X)\eta(Y)+2(m-1)\eta(X)g(\phi SZ,Y)\nonumber\\
&-2g(SZ,AX)g(AY,N)+4q(Z)g(N,AX)g(AY,N)\nonumber\\
&-2g(SZ,AY)g(AX,N)-g(SZ, S\phi SX)\eta(Y)+g(\phi(\nabla_ZS)\phi SX,Y)\\
&+\eta(SX)g(\phi S^2Z,Y)+g(\phi S\phi(\nabla_ZS)X,Y).\nonumber
\end{align}
Since $S\xi=\alpha\xi$, letting  $X=\xi$ we get
\begin{align*}
  0=&2(m-1)g(\phi SZ,Y)-2g(SZ,A\xi)g(AY,N)\\
&+\alpha g(\phi S^2Z,Y)+g((\nabla_ZS)\xi,\phi S\phi Y)\\
=&2(m-1)g(\phi SZ,Y)-2g(SZ,A\xi)g(AY,N)\\
&+\alpha g(\phi S^2Z,Y)+g(\alpha\phi SZ-S\phi SZ,\phi S\phi Y).
\end{align*}
Moreover, if $Z=\xi$ then we get $\alpha g(A\xi,\xi)g(AY,N)=0$.
If $\alpha\neq0$ then $\cos(2t)g(AY,N)=0$ by \eqref{9}.
 That means that $t=\frac{\pi}{4}$ or $AY\in TM$, that is,
 the unit normal vector $N$ is $\mathfrak{U}$-principal
 or $\mathfrak{U}$-isotropic. If $\alpha=0$ then $g(Y,AN)g(\xi,A\xi)=0$
 for any $Y\in TM$ by \eqref{14}, thus we have same conclusion. The proof is complete.\qed\bigskip

We first assume that the unit normal vector field $N$ is $\mathfrak{U}$-isotropic.
In this case these expressions in \eqref{9} become
\begin{equation*}
\left\{
\begin{array}{ll}
N&=\frac{1}{\sqrt{2}}(Z_1+JZ_2),\\
AN&=\frac{1}{\sqrt{2}}(Z_1-JZ_2),\\
\xi&=\frac{1}{\sqrt{2}}(Z_2-JZ_1),\\
A\xi&=\frac{1}{\sqrt{2}}(Z_2+JZ_1).
  \end{array}
\right.
\end{equation*}
Thus
\begin{equation}\label{22}
  g(A\xi,\xi)=g(AN,N)=0.
\end{equation}
So \eqref{15} becomes
\begin{align}\label{23}
S\phi SX=&\frac{1}{2}\alpha (\phi S+S\phi )X+\phi X\nonumber\\
&-g(X,AN)A\xi+g(X,A\xi)AN.
 \end{align}
The formula \eqref{21} with $Z=\xi$ implies
\begin{align}\label{24}
  0=&-2g(S\xi,AX)g(AY,N)+4q(\xi)g(N,AX)g(AY,N)\nonumber\\
&-2g(S\xi,AY)g(AX,N)-g((\nabla_\xi S)\phi SX,\phi Y)\\
&+g((\nabla_\xi S)X,\phi S\phi Y).\nonumber
\end{align}
By Codazzi equation \eqref{8}, we get
\begin{align*}
 (\nabla_\xi S)Y=&\alpha\phi SY-S\phi SY+\phi Y-g(Y,AN)A\xi\\
&+g(Y,A\xi)AN\\
=&\frac{1}{2}\alpha(\phi S-S\phi)Y.
 \end{align*}
Thus substituting this into \eqref{24} gives
\begin{equation}\label{25}
\begin{aligned}
  0=&-2\alpha g(\xi,AX)g(AY,N)+4q(\xi)g(N,AX)g(AY,N)\\
&-2\alpha g(\xi,AY)g(AX,N)-\frac{1}{2}\alpha g(S\phi SX+\phi S\phi S\phi X, Y).
\end{aligned}\end{equation}
Moreover, by \eqref{23} we have $S\phi SX+\phi S\phi S\phi X=0$, thus taking $X=A\xi$ in \eqref{25} yields
\begin{align*}
  \alpha g(AY,N)=0.
\end{align*}
Here we have used $g(A\xi,A\xi)=1$ and $g(AN,A\xi)=0.$
From this we derive $\alpha=0$ since $N$ is $\mathfrak{U}$-isotropic.

On the other hand, we put $Y=\xi$ in \eqref{21} and get
\begin{align*}
  0=&2(m-1)g(SZ,\phi X)+2g(SZ,A\xi)g(AX,N)+g(SZ, S\phi SX).
\end{align*}
 Applying \eqref{23} in the above formula, we have
\begin{align*}
0=&(2m-1)g(SZ,\phi X)+g(SZ,A\xi)g(AX,N)+g(SZ,AN)g(X,A\xi).
 \end{align*}
 That is,
\begin{align}\label{26}
0=&(2m-1)S\phi X+g(AX,N)SA\xi+g(X,A\xi)SAN.
 \end{align}
When $X=AN$, it comes to
\begin{align*}
0=&(2m-1)S\phi AN+SA\xi.
 \end{align*}
Then $A\xi=\phi AN$ implies $SA\xi=0$. Similarly, $SAN=0.$
Therefore from \eqref{26} we obtain $S\phi X=0$ for all $X\in TM.$
As $S\xi=0$ we know $SX=0$ for all $X\in TM$, thus $\nabla_\xi S=0$,
that means that the hypersurface $M$ admits  parallel shape operator.
 But Suh \cite{Suh2} has showed the non-existence of this type hypersurfaces.

In the following if $N$ is $\mathfrak{U}$-principal, that is, $AN=N$, then \eqref{13} becomes
\begin{align*}
  \mathrm{Ric}^*(X,Y) =&-2(m-1)g(\phi^2 X,Y)-g((\phi S)^2X,Y).
\end{align*}
In this case we see that the star-Ricci tensor is commuting by Lemma \ref{L2.1}.
Thus we see $\alpha=0$ from the proof of Theorem \ref{T1*}.
In this case, the formulas \eqref{21} with $Y=\xi$ and \eqref{15} respectively become
$2(m-1)g(SZ,\phi X)+g(SZ,S\phi SX)=0$ and $S\phi SX=\phi X$, respectively.
From these two equations we obtain $g(SZ,\phi X)=0,$
that is, $\phi SZ=0.$ This implies $SZ=\alpha\eta(Z)\xi=0.$ As before, this is impossible.

Summing up the above discussion, we complete the proof Theorem \ref{T3}.

\section{Proof of Theorem \ref{T2}}
In order to prove our theorem, we first give the following property.
\begin{prop}\label{P1}
Let $M$ be a real hypersurface in $Q^m,m\geq3$, admitting a star-Ricci
soliton with potential vector field $\xi$, then $M$ must be Hopf.
\end{prop}
\proof Since $\mathcal{L}_Wg$ and $g$ are symmetry, the *-Ricci soliton
equation \eqref{2} implies the star-Ricci tensor is also symmetry,
i.e. $\mathrm{Ric}^*(X,Y)=\mathrm{Ric}^*(Y,X)$ for any vector fields $X,Y$ on $M$. It yields from \eqref{13} that
\begin{equation}\label{27}
(\phi S)^2X=(S\phi)^2X
\end{equation}
for all $X\in TM.$

On the other hand, from the star-Ricci soliton equation \eqref{2} it follows
\begin{equation}\label{3.7}
  \mathrm{Ric}^*(X,Y)=\lambda g(X,Y)+\frac{1}{2}g((S\phi-\phi S)X,Y).
\end{equation}
By \eqref{13}, we have
\begin{align}\label{29}
  &-2(m-1)g(\phi^2 X,Y)-2g(N,AX)g(AY,N)-g((\phi S)^2X,Y)\nonumber\\
  &=\lambda g(X,Y)+\frac{1}{2}g((S\phi-\phi S)X,Y).
\end{align}
Putting $X=Y=\xi$ gives $\lambda=0$ since $g(AN,\xi)=0.$
Therefore the previous formula with $X=\xi$ yields
\begin{align*}
  (\phi S)^2\xi=\frac{1}{2}\phi S\xi.
\end{align*}
Using \eqref{27} we get $\phi S\xi=0$, which shows
$S\xi=\alpha\xi$ with $\alpha=g(S\xi,\xi).$\qed\bigskip

Moreover, by \eqref{3.7} we have
 \begin{align}\label{30}
\mathrm{Ric}^*(X)=\frac{1}{2}(S\phi-\phi S)X.
\end{align}
Thus by a straightforward computation we find $\phi\circ \mathrm{Ric}^*+\mathrm{Ric}^*\circ\phi=0$
since the relation $\phi^2S=S\phi^2$ holds by Proposition \ref{P1}. Namely the following result holds.
\begin{prop}
Let $M$ be a real hypersurface in $Q^m,m\geq3$, admitting a
star-Ricci soliton with potential vector field $\xi$, then the star-Ricci tensor is anti-commuting.
\end{prop}

Next we will compute the convariant derivative of $\phi\circ \mathrm{Ric}^*+\mathrm{Ric}^*\circ\phi=0.$
First of all, by \eqref{30} and \eqref{6}, we compute
  \begin{align}\label{31}
(\nabla_X\mathrm{Ric}^*)(Y)
=&\frac{1}{2}\Big\{(\nabla_XS)\phi Y+S(\nabla_X\phi) Y-(\nabla_X\phi) SY-\phi(\nabla_XS)Y\Big\}\nonumber\\
=&\frac{1}{2}\Big\{(\nabla_XS)\phi Y+\eta(Y)S^2X-\alpha g(SX,Y)\xi\nonumber\\
&-\alpha\eta(Y)SX+g(SX,SY)\xi-\phi(\nabla_XS)Y\Big\}.
\end{align}
Now differentiating $\phi\circ\mathrm{Ric}^*+\mathrm{Ric}^*\circ\phi=0$ convariantly gives
\begin{align*}
   0=&(\nabla_X\phi) \mathrm{Ric}^*(Y)+\phi(\nabla_X\mathrm{Ric}^*)Y+(\nabla_X\mathrm{Ric}^*)\phi Y+\mathrm{Ric}^*(\nabla_X\phi)Y  \\
    =& -g(SX,\mathrm{Ric}^*(Y))\xi+\phi(\nabla_X\mathrm{Ric}^*)Y+(\nabla_X\mathrm{Ric}^*)\phi Y+\eta(Y)\mathrm{Ric}^*(SX)\\
=& -\frac{1}{2}g(SX,S\phi Y-\phi SY))\xi+\phi(\nabla_X\mathrm{Ric}^*)Y+(\nabla_X\mathrm{Ric}^*)\phi Y\\
&+\frac{1}{2}\eta(Y)(S\phi SX-\phi S^2X).
 \end{align*}
Applying \eqref{31} in the above formula, we get
\begin{align*}
0=& g(SX,\phi SY)\xi+\Big\{-\alpha\eta(Y)\phi SX+g((\nabla_XS)Y,\xi)\xi\Big\}\\
&+\Big\{\eta(Y)(\nabla_XS)\xi-\alpha g(SX,\phi Y)\xi\Big\}+\eta(Y)S\phi SX\\
=& g(SX,\phi SY)\xi-\alpha\eta(Y)\phi SX+\Big\{g((Y,X(\alpha)\xi+\alpha\phi SX-S\phi SX)\Big\}\xi\\
&+\eta(Y)\Big\{X(\alpha)\xi+\alpha\phi SX-S\phi SX\Big\}-\alpha g(SX,\phi Y)\xi+\eta(Y)S\phi SX\\
=& 2g(SX,\phi SY)\xi+2\eta(Y)X(\alpha)\xi-2\alpha g(SX,\phi Y)\xi,
 \end{align*}
i.e.
\begin{equation}\label{32}
g(SX,\phi SY)+\eta(Y)X(\alpha)-\alpha g(SX,\phi Y)=0.
\end{equation}
 From this we know $X(\alpha)=0$ by taking $Y=\xi,$
i.e. $\alpha$ is constant. Hence formula \eqref{32} becomes
\begin{equation*}
g(SX,\phi SY)=\alpha g(SX,\phi Y).
\end{equation*}
Now interchanging $X$ and $Y$  and comparing the resulting equation with
the previous equation, we have $\alpha(\phi S-\phi S)X=0$, which shows that either $\alpha=0$ or $\phi S=S\phi.$
 Namely the following lemma has been proved.
\begin{lemma}
Let $M$ be a real hypersurface in $Q^m,m\geq3$, admitting a star-Ricci
soliton with potential vector field $\xi$, then either the Reeb flow is isometric, or $\alpha=0$.
\end{lemma}

If the Reeb flow of $M$ is isometric, Berndt and Suh proved the following conclusion:
\begin{thm}[\cite{BS2}]\label{B}
 Let $M$ be a real hypersurface of the complex quadric $Q^m,m\geq3.$ The
Reeb flow on $M$ is isometric if and only if $m$ is even, say $m=2k$, and $M$ is an open part
of a tube around a totally geodesic $\mathbb{C}P^k\subset Q^{2k}$.
\end{thm}

In the following we set $\alpha=0$, it follows from \eqref{32} that
\begin{equation}\label{33}
  S\phi SX=0,\quad\text{for all}\quad X\in TM.
\end{equation}
And it is easy to show that the normal vector $N$ is either
$\mathfrak{U}$-principal or $\mathfrak{U}$-isotropic from \eqref{14}.
In the following let us consider these two cases.

{\bf Case I:} $N$ is $\mathfrak{U}$-principal, that is, $AN=N.$ We follow from \eqref{15} that
\begin{equation*}
  S\phi SX=\phi X.
\end{equation*}
By comparing with \eqref{33} we find $\phi X=0$, which is impossible.

{\bf Case II:} $N$ is $\mathfrak{U}$-isotropic.
Using \eqref{33}, we derive from \eqref{15}
\begin{equation}\label{34}
g(\phi X,Y)=g(X,AN)g(Y,A\xi)-g(Y,AN)g(X,A\xi).
\end{equation}
Using \eqref{33} again, we learn \eqref{29} becomes
\begin{align*}
  -&2(m-1)g(\phi^2 X,Y)-2g(N,AX)g(AY,N)\nonumber\\
  &=\frac{1}{2}g((S\phi-\phi S)X,Y).
\end{align*}
Moreover, replacing $Y$ by $\phi Y$ gives
\begin{align}\label{35}
  -&2(m-1)g(\phi X,Y)+2g(N,AX)g(Y,A\xi)\nonumber\\
  &=\frac{1}{2}g((S\phi-\phi S)X,\phi Y).
\end{align}
Here we have used $g(A\phi Y,N)=-g(Y,A\xi)$, which follows from \eqref{19} and \eqref{22}.

By interchanging $Y$ and $X$ in the formula \eqref{35} and
applying the resulting equation to subtract this equation, we get
 \begin{align*}
   &2g(N,AX)g(Y,A\xi)-2g(N,AY)g(X,A\xi)  \\
=&\frac{1}{2}g((S\phi-\phi S)X,\phi Y)+2(m-1)g(\phi X,Y) \\
    & -\frac{1}{2}g((S\phi-\phi S)Y,\phi X)-2(m-1)g(\phi Y,X)\\
=&4(m-1)g(\phi X,Y).
 \end{align*}
Combining this with \eqref{34} we get $(m-3)\phi X=0,$ which is a
contradiction if $m\geq4.$ Hence we complete the proof of Theorem \ref{T2}.\bigskip\\
{\bf Acknowledgement.}
The author would like to thank the referee for the
valuable comments on this paper.


\begin{thebibliography}{9}

\bibitem{BS} \textsc{J. Berndt, Y. J. Suh},  Hypersurfaces in Kaehler manifold, Proc. A.M.S. {\bf143} (2015), 2637-2649.
\bibitem{BS2} \textsc{J. Berndt, Y. J. Suh},  Real hypersurfaces with isometric Reeb flows in complex quadrics,
Inter. J. Math. {\bf24} (2013), 1350050, 18pp.

\bibitem{HA}\textsc{T. Hamada},  Real hypersurfaces of complex space forms
in terms of Ricci *-tensor, Tokyo J. Math. {\bf 25} (2002), 473-483.

\bibitem{IR}\textsc{T. A. Ivey, P. J. Ryan}, The $^{*}$-Ricci tensor
for Hypersurfaces in $\mathbb{CP}^n$ and $\mathbb{C
}\mathrm{H}^n$, Tohoku Math. J. {\bf34} (2011), 445-471.

\bibitem{KP} \textsc{G. Kaimakamis, K. Panagiotidou},  *-Ricci solitons
of real hypersurfaces in non-flat complex space forms,
            J. Geom. Phys. {\bf 86} (2014), 408--413.

\bibitem{K} \textsc{S. Klein}, Totally geodesic submanifolds
in the complex quadric, Diff. Geom. Appl. {\bf26} (2008) 79-96.

\bibitem{KN} \textsc{S. Kobayashi, K. Nomizu},  Foundations of Differential Geometry,
Vol. II, Wiley Classics Library ed., A Wiley-Interscience Publ., 1996.


\bibitem{R} \textsc{H. Reckziegel},  On the geometry of the complex quadric, in:
Geometry and Topology of Submanifolds VIII, Brussels/Nordfjordeid, 1995,
World Sci. Publ., River Edge, NJ, 1995, 302-315.

\bibitem{Suh1} \textsc{Y. J. Suh},  Real hypersurfaces in the complex
quadric with parallel Ricci tensor, Adv. Math. {\bf281} (2015), 886-905.

\bibitem{Suh2} \textsc{Y. J. Suh},  Real hypersurfaces in the complex
quadric with Reeb parallel shape operator, Inter. J. Math. {\bf 25} (2014) 1450059, 17pp.

\bibitem{Suh3} \textsc{Y. J. Suh},  Real hypersurfaces in the complex quadric with parallel normal
Jacobi operator, Math. Nachr. {\bf289} (2016), 1-10.

\bibitem{Suh4} \textsc{Y. J. Suh},  Real hypersurfaces in the complex quadric
with harmonic curvature, J. Math. Pures. Appl. {\bf106} (2016), 393-410.

\bibitem{Suh5}  \textsc{Y. J. Suh},  Real hypersurfaces in the complex quadric with commuting
and parallel Ricci tensor, J.  Geom. Phys. {\bf106} (2016), 130-142.

\bibitem{Suh6} \textsc{Y. J. Suh},  Pseudo-anti commuting Ricci
tensor and Ricci soliton real hypersurfaces in the
complex quadric, J. Math. Pure. Appl. {\bf107} (2017), 429-450.

\bibitem{T}\textsc{S. Tachibana}, On almost-analytic vectors
in almost-K\"ahlerian manifolds, Tohoku Math. J. {\bf11} (1959), 247-265.

\end{thebibliography}
\end{document}